\documentclass{amsart}
\date{September 3, 2012}
\usepackage{hyperref}
\usepackage{latexsym,amsmath,amsfonts,amscd,amssymb}

\theoremstyle{plain}  
\newtheorem{theorem}{Theorem}[section]

\newtheorem*{theorem*}{Theorem}

\newtheorem{proposition}[theorem]{Proposition}

\theoremstyle{definition}
\newtheorem{definition}[theorem]{Definition}

\theoremstyle{remark}
\newtheorem{example}[theorem]{Example}

\newtheorem*{notation*}{Notation}
\newtheorem{remark}[theorem]{Remark}

\newtheorem*{question*}{Question}

\newtheorem*{claim*}{Claim}

\numberwithin{equation}{section}

\newcommand{\suchthat}{\;\;|\;\;}
\newcommand{\abs}[1]{\lvert#1\rvert}
\newcommand{\norm}[1]{\lVert#1\rVert}

\renewcommand{\leq}{\leqslant}

\renewcommand{\geq}{\geqslant}

\newcommand{\into}{\hookrightarrow}

\newcommand{\xra}{\xrightarrow}

\newcommand{\lie}{\mathfrak}

\newcommand{\R}{\mathbb{R}}

\newcommand{\Z}{\mathbb{Z}}
\newcommand{\C}{\mathbb{C}}

\newcommand{\HH}{\mathbb{H}}

\newcommand{\dbar}{\bar{\partial}}

\newcommand{\SU}{\mathrm{SU}}
\newcommand{\U}{\mathrm{U}}

\newcommand{\GL}{\mathrm{GL}}
\newcommand{\SL}{\mathrm{SL}}

\newcommand{\SO}{\mathrm{SO}}
\newcommand{\Sp}{\mathrm{Sp}}

\newcommand{\Lie}{\mathrm{Lie}}
\newcommand{\x}{\times}

 \DeclareMathOperator{\ad}{ad}
\DeclareMathOperator{\Ad}{Ad} 
\DeclareMathOperator{\grad}{grad} 
\DeclareMathOperator{\rk}{rk} 
 
\DeclareMathOperator{\Hom}{Hom} \DeclareMathOperator{\End}{End}

\DeclareMathOperator{\Index}{index}

\newcommand{\aut}{\operatorname{aut}}
\newcommand{\Aut}{\operatorname{Aut}}

\newcommand{\MB}{\mathcal{M}^{\mathrm{B}}}
\newcommand{\MdR}{\mathcal{M}^{\mathrm{dR}}}
\newcommand{\MDol}{\mathcal{M}^{\mathrm{Dol}}}
\newcommand{\Mhar}{\mathcal{M}^{\mathrm{Har}}}
\newcommand{\MHit}{\mathcal{M}^{\mathrm{Hit}}}
\newcommand{\wt}{\widetilde}
\renewcommand{\phi}{\varphi}

\newcommand{\liem}{\mathfrak{m}}

\newcommand{\liemc}{\mathfrak{m}^{\mathbb{C}}}
\newcommand{\lieh}{\mathfrak{h}}
\newcommand{\liehc}{\mathfrak{h}^{\mathbb{C}}}
\newcommand{\lieg}{\mathfrak{g}}
\newcommand{\liegc}{\mathfrak{g}^{\mathbb{C}}}
\newcommand{\vol}{\mathrm{vol}}

\newcommand{\CC}{\mathbb{C}}

\newcommand{\QQ}{\mathbb{Q}}

\begin{document}

\title{Surface group representations and Higgs bundles}

\author[P. B. Gothen]{Peter B. Gothen}
\address{Centro de Matem\'atica \\
Faculdade de Ci\^encias da Universidade do Porto \\
Rua do Campo Alegre \\ 4169-007 Porto \\ Portugal }
\email{pbgothen@fc.up.pt}


\thanks{
Member of VBAC (Vector Bundles on Algebraic Curves).
Partially supported by the FCT (Portugal) with EU (COMPETE) and national funds 
through the projects PTDC/MAT/099275/2008 and PTDC/MAT/098770/2008,
and through Centro de Matem\'atica da Universidade do Porto
(PEst-C/MAT/UI0144/2011).
}


\maketitle


\section{Introduction}
\label{sec:introduction}

In these notes we give an introduction to Higgs bundles and their
application to the study of surface group representations. This is
based on two fundamental theorems. The first is the theorem of
Corlette and Donaldson on the existence of harmonic metrics in flat
bundles which we treat in Lecture~1, after explaining some
preliminaries on surface group representations, character varieties and
flat bundles. The second is the Hitchin--Kobayashi correspondence for
Higgs bundles, which goes back to the work of Hitchin and Simpson;
this is the main topic of Lecture 2. Together, these two results
allows the character variety for representations of the fundamental
group of a Riemann surface in a Lie group $G$ to be identified with
a moduli space of holomorphic objects, known as $G$-Higgs
bundles. Finally, in Lecture~3, we show how the $\C^*$-action on the
moduli space $G$-Higgs bundles can be used to study its topological
properties, thus giving information about the corresponding character
variety.

For lack of time and expertise, we do not treat many other important
aspects of the theory of surface group representations, such as the
approach using bounded cohomology (see, e.g., Burger--Iozzi--Wienhard
\cite{burger-iozzi-wienhard:2010b,burger-iozzi-wienhard:2010}), higher
Teichm\"uller theory (see, e.g., Fock--Goncharov
\cite{fock-goncharov:2006}), or ideas related to geometric structures
on surfaces (see, e.g., Goldman \cite{goldman:2010}). We also do not
touch on the relation of Higgs bundle moduli with mirror symmetry and
the Geometric Langlands Programme (see, e.g., Hausel
\cite{hausel:2011} and Kapustin--Witten \cite{kapustin-witten:2007}).

In keeping with the lectures we do not give proofs of most
results. For more details and full proofs, we refer to the
literature. Some references that the reader may find useful are the
papers of Hitchin \cite{hitchin:1987a,hitchin:1992},
Garc\'\i{}a-Prada \cite{garcia-prada:2009}, Goldman
\cite{goldman:1985,goldman:2010,goldman:2010b} and also the papers
\cite{bradlow-garcia-prada-gothen:2003,bradlow-garcia-prada-gothen:hss-higgs,garcia-gothen-mundet:2009b,garcia-gothen-mundet:2009a,bradlow-garcia-gothen:2009}.

\subsection*{Notation}

Our notation is mostly standard. Smooth $p$-forms are denoted by
$\Omega^p$ and smooth $(p,q)$-forms by $\Omega^{p,q}$. We shall
occasionally confuse vector bundles and locally free sheaves.

\subsection*{Acknowledgments}

It is a pleasure to thank the organizers of the Newton Institute
Semester on Moduli Spaces, and especially Leticia Brambila-Paz, for
the invitation to lecture in the School on Moduli spaces and for
making it such a pleasant and stimulating event. I would also like to
thank the participants in the course for their interest and for making
the tutorials a fun and rewarding experience.
It is impossible to mention all the mathematicians to whom I am
indebted and who have generously shared their insights on the topics
of these lectures over the years and without whom these notes would
not exist. But I would like to express my special gratitude
to Bill Goldman, Nigel Hitchin and my collaborators Ignasi Mundet i
Riera, Steve Bradlow and Oscar Garc\'\i{}a-Prada.

\section{Lecture 1: Character varieties
  for surface groups and harmonic maps}

In this lecture we give some basic definitions and properties of
character varieties for representations of surface groups. We then
explain the theorem of Corlette and Donaldson on the existence of
harmonic maps in flat bundles, which is one of the two central results
in the non-abelian Hodge theory correspondence (the other one being
the Hitchin--Kobayashi correspondence, which will be treated in
Lecture 2). 
 
\subsection{Surface group representations and character varieties}
More details on the following can be found in, for example,
Goldman\cite{goldman:1985}.

Let $\Sigma$ be a closed oriented surface of genus $g$. The
fundamental group of $\Sigma$ has the standard presentation
\begin{equation}
  \label{eq:1}
  \pi_1 \Sigma = \langle a_1,b_1,\dots,a_g,b_g \suchthat
  \prod[a_i,b_i]=1\rangle,
\end{equation}
where $[a_i,b_i] = a_ib_ia_i^{-1}b_i^{-1}$ is the commutator.

Let $G$ be a real reductive Lie group. We denote its Lie algebra by
$\lie{g} = \Lie(G)$. Though not strictly necessary for everything that
follows, we shall assume that $G$ is connected. We shall also fix a
non-degenerate quadratic form on $G$, invariant under the adjoint
action of $G$ (when $G$ is semisimple, the Killing form or a multiple
thereof will do).

By definition a \emph{representation} of $\pi_1 \Sigma$ in $G$ is a
homomorphism $\rho\colon \pi_1 \Sigma \to G$. Let $\Ad\colon G \to
\Aut(\lieg)$ be the adjoint representation of $G$ on its Lie algebra
$\lieg$. We say that $\rho$ is \emph{reductive}\footnote{When $G$ is
  algebraic an alternative equivalent definition is to ask for the
  Zariski closure of $\rho(\Sigma) \subset G$ to be reductive.} if the
composition
\begin{displaymath}
  \Ad \circ \rho\colon \pi_1\Sigma \to \Aut(\lieg)
\end{displaymath}
is completely reducible. Denote by $\Hom^{\mathrm{red}}(\pi_1\Sigma,G)
\subset \Hom(\pi_1\Sigma,G)$ the subset of reductive representations.

\begin{definition}
The \emph{character variety} for representations of $\pi_1 \Sigma$ in
$G$ is
\begin{displaymath}
  \MB(\Sigma,G) = \Hom^{\mathrm{red}}(\pi_1\Sigma,G) / G,
\end{displaymath}
where the $G$-action is by simultaneous conjugation:
\begin{displaymath}
  (g\cdot\rho)(x) = g\rho(x)g^{-1}.
\end{displaymath}
\end{definition}

The character variety is also known as the \emph{Betti
  moduli space} (in Simpson's language \cite{simpson:1992}).

Note that, using the presentation (\ref{eq:1}), a representation $\rho$
is given by a $2g$-tuple of elements in $G$ satisfying the
relation. Hence we get an inclusion $\Hom(\pi_1\Sigma,G) \into G^{2g}$, which
endows $\Hom(\pi_1\Sigma,G)$ with a natural topology.\footnote{This is
  in fact the same as the compact-open topology on the mapping space
  $\Hom(\pi_1\Sigma,G)$, where we give $\pi_1\Sigma$ the discrete
  topology.}  However, it turns out that the quotient space $\Hom(\pi_1,G) / G$
is not in general Hausdorff. The restriction to reductive
representations remedies this problem.

We also remark that, in case $G$ is a complex reductive algebraic
group, the character variety can be constructed as an affine GIT
quotient (this is classical; a nice exposition is contained in \S 3.1
of Casimiro--Florentino \cite{casimiro-florentino:2011}).

\subsection{Review of connections and curvature in principal bundles}

Recall that a (smooth) \emph{principal $G$-bundle} on $\Sigma$ is a
smooth fibre bundle $\pi\colon E \to \Sigma$ with a $G$-action (normally taken
to be on the right) which is free and transitive on each
fibre. Moreover, $E$ is required to 
admit $G$-equivariant local trivializations $E_{|U} \cong U \times
G$ over small open sets $U \subset \Sigma$ (where $G$ acts by right
multiplication on the second factor of the product $U \times G$). Note
that the fibre $E_x$ over any $x \in \Sigma$ is a $G$-torsor so,
choosing an element $e \in E_x$, we get a canonical identification
$E_x \cong G$.

\begin{example}
  (1) The \emph{frame bundle} of a rank $n$ complex vector bundle $V\to
  \Sigma$ is a principal $\GL(n,\C)$-bundle, which has fibres
  \begin{displaymath}
    E_x=\{e\colon \C^n \xra{\cong} V_x \suchthat \text{$e$ is a linear
    isomorphism}\}.
  \end{displaymath}
 
  (2) The universal covering $\widetilde\Sigma \to \Sigma$ is a
  principal $\pi_1\Sigma$-bundle over $\Sigma$. In this case
  the action is on the left.
\end{example}

Whenever we have a principal $G$-bundle $E\to \Sigma$ and a smooth
$G$-space $V$ (i.e., $V$ is a smooth manifold on which $G$ acts by
smooth maps), we obtain a fibre bundle $E(V)$
with fibres modeled on $V$ by taking the quotient of $E\times V$
under the diagonal $G$-action:
\begin{displaymath}
  E(V) = E\times_G V \to \Sigma.
\end{displaymath}
In particular, if $V$ is a vector
space $V$ with a linear $G$-action, we obtain a vector bundle $E(V)$
with fibres modeled on $V$.
An important instance of this construction is when $V = \lieg$ acted
on by $G$ via the adjoint action. The resulting vector bundle $\Ad E :=
E(\lieg)$ is then known as the \emph{adjoint bundle} of $E$.

There is a bijective correspondence between sections $s\colon\Sigma\to
E(V)$ of
the bundle $\pi\colon E(V)\to \Sigma$ and $G$-equivariant maps $\tilde{s}\colon
E \to V$, given by
\begin{displaymath}
  s(x) = [e,\tilde{s}(e)]
\end{displaymath}
for $e \in E(V)_x = \pi^{-1}(x)$ and $x \in \Sigma$.
Similarly, a $G$-equivariant differential $p$-form
$\alpha\in\Omega^p(E,V)$ descends to an $E(V)$-valued $p$-form
$\tilde\alpha\in\Omega^p(\Sigma,E(V))$ if and only if it is \emph{tensorial},
i.e., it vanishes on the vertical tangent spaces $T^v_eE = T_eE_x$ to
$E$.

A \emph{connection} in a principal $G$-bundle $E\to \Sigma$ is given
by a smooth $G$-invariant Lie algebra valued $1$-form $A \in
\Omega^1(E,\lieg)$ which restricts to the identity on the vertical
tangent spaces $T^v_eE$ under the natural identification $T^v_eE \cong
\lieg$ given by the choice of $e \in E_x$. Equivalently, a connection
corresponds to the choice of a horizontal complement $T^h_eE =
\ker(A(e)\colon T_eE \to \lieg)$ to $T^v_eE$ in each $T_eE$. Moreover,
the $G$-invariance means that these complements correspond under the
$G$-action. The difference of two connections is a tensorial form, so
it follows that the space $\mathcal{A}$ of connections on $E$ is an
affine space modeled on $\Omega^1(\Sigma, \Ad E)$.

Given a connection $A$ in a principal bundle $E$, we obtain a
covariant derivative
\begin{displaymath}
  d_A\colon \Omega^0(\Sigma,E(V)) \to \Omega^1(\Sigma,E(V))
\end{displaymath}
on sections in any associated vector bundle $E(V)$ as follows. 
Let $s\in \Omega^0(\Sigma,E(V))$ and let $\tilde{s}\colon E \to V$ be
the corresponding $G$-equivariant map as above. Then we define a
tensorial one-form $\wt{d_A(s)}$ on $E$ by composing $d\tilde{s}$ with
the projection $TE \to T^hE$ defined by $A$, and let
$d_A(s)\in\Omega^1(\Sigma,E(V))$ be the corresponding $E(V)$-valued
one-form. 

Given a connection in $E$, the horizontal subspaces define a
$G$-invariant distribution on the total space of $E$. The
obstruction to integrability of this horizontal distribution is given
by the \emph{curvature}
\begin{displaymath}
  F(A) = dA + \tfrac{1}{2}[A,A] \in \Omega^2(E,\lieg)
\end{displaymath}
of the connection $A$, where the bracket $[A,A]$ is defined by
combining the wedge product on forms with the Lie bracket on $\lieg$.
One checks that $F(A)$ is in fact a tensorial form and therefore
descends to a $2$-form on $\Sigma$ with values in the adjoint bundle,
which we denote by the same symbol,
\begin{displaymath}
  F(A) \in \Omega^2(\Sigma,E(\lieg)).
\end{displaymath}
A connection $A$ is \emph{flat} if $F(A) = 0$. A principal $G$-bundle
$E \to \Sigma$ with a flat connection is called a \emph{flat
  bundle}. Equivalently, a flat bundle is one for which the structure
group $G$ is discrete. 
The Frobenius Theorem has the following immediate consequence.

\begin{proposition}
  \label{prop:frob}
  Let $E \to \Sigma$ be a flat bundle and let $e_0 \in E_{x_0}$ for
  some $x_0 \in \Sigma$. Then, for any sufficiently small neighbourhood
  $U \subset \Sigma$, there is a unique section $s \in
  \Omega^0(U,E_{|U})$ such that $d_A(s) = 0$ and $s(x_0) = e_0$.
\end{proposition}

\subsection{Surface group representations and flat bundles}

Given a $G$-bundle $E$ on $\Sigma$ with a connection $A$, it follows
from the existence and uniqueness theorem for ordinary differential
equations that we can lift any loop $\gamma$ in $\Sigma$ to a
covariantly constant loop in $E$ (i.e., one whose tangent vectors are
horizontal for the connection). In this way we obtain a well-defined
parallel transport $E_x\to E_x$, which is given by multiplication by a
unique group element, the \emph{holonomy of $A$ along $\gamma$},
denoted by $h_A(\gamma)\in G$. Moreover, if the connection $A$ is
flat, it follows from Proposition~\ref{prop:frob} that the holonomy
only depends on the homotopy class of $\gamma$ and thus we obtain the
\emph{holonomy representation} of $\pi_1\Sigma$:
\begin{equation}
  \rho_A\colon\pi_1\Sigma \to G
\end{equation}
defined by $\rho_A([\gamma]) = h_A(\gamma)$. We say that a flat
connection $A$ is \emph{reductive} if its holonomy representation is a
reductive representation of $\pi_1\Sigma$ in $G.$
 
On the other hand, let $\rho\colon \pi_1\Sigma \to G$ be
representation. We can then define a principal $G$-bundle $E_\rho$ by
taking the quotient
\begin{displaymath}
  E_\rho = \widetilde\Sigma \times_{\pi_1\Sigma}G,
\end{displaymath}
where $\pi_1\Sigma$ acts on the universal cover $\widetilde \Sigma \to
\Sigma$ by deck
transformations and on $G$ by left multiplication via
$\rho$. Moreover, since $\widetilde \Sigma \to \Sigma$ is a covering,
there is a natural choice of horizontal subspaces in
$E_\rho$. Therefore this bundle has a naturally defined connection
which is evidently flat. 

One sees that these two constructions are inverses of each other. Next
we shall introduce the natural equivalence relation on (flat)
connections and promote this correspondence to a bijection between
equivalence classes of flat connections and points in the character
variety.

\subsection{Flat bundles and gauge equivalence}

The \emph{gauge group}\footnote{This is the mathematician's
  definition. To a physicist the gauge group is the structure group $G$.} 
is the automorphism group
\begin{displaymath}
  \mathcal{G} = \Omega^0(\Sigma,\Aut(E))
\end{displaymath}
where $\Aut(E) = E\x_{\Ad}G \to \Sigma$ is the bundle of automorphisms
of $E$. The gauge group acts on the space of connections
$\mathcal{A}_E$ via
\begin{displaymath}
  g\cdot A = gAg^{-1}+gdg^{-1}.
\end{displaymath}
Moreover, the corresponding action on the curvature is
\begin{equation}
  \label{eq:2}
  F(g\cdot A) = gF(A)g^{-1}
\end{equation}
and hence $\mathcal{G}$ preserves the subspace of flat connections
on $E$. 

Recall that principal $G$-bundles on $\Sigma$ are classified (up to
smooth isomorphism) by a
characteristic class
\begin{displaymath}
  c(E) \in H^2(\Sigma,\pi_1 G) \cong \pi_1 G.
\end{displaymath}
Here we are using the fact that $G$ is connected and that $\Sigma$ is
a closed oriented surface.  Fix $d \in \pi_1 G$ and let $E\to \Sigma$
be a principal $G$-bundle with $c(E)=d$. We can then consider the
quotient space
\begin{displaymath}
  \MdR_d(\Sigma,G) = \{ A \in \mathcal{A} \suchthat \text{$F(A) = 0$ and
  $A$ is reductive}\} / \mathcal{G},
\end{displaymath}
which is known as the \emph{de Rham moduli space} (recall that
$\mathcal{A}$ denotes the space of connections).


\begin{proposition}
  \label{prop:flat-rep}
  If flat connections $B_i$ correspond to representations
  $\rho_i\colon\pi_1\Sigma\to G$ for $i=1,2$, then $B_1$ and $B_2$ are
  gauge equivalent if and only if there is a $g \in G$ such that
  $\rho_1 = g\rho_2g^{-1}$.
\end{proposition}

This proposition implies that there is a bijection
\begin{equation}
  \label{eq:21}
  \MdR_d(\Sigma,G) \cong \MB_d(\Sigma,G),
\end{equation}
where we denote by $\MB_d(\Sigma,G) \subset \MB(\Sigma,G)$ the subspace
of representations with characteristic class $d$.




\subsection{Harmonic metrics in flat bundles}

Let $G' \subset G$ be a Lie subgroup. Recall that a \emph{reduction of
  structure group} in a principal $G$-bundle $E \to \Sigma$ to $G'
\subset G$ is a section
\begin{displaymath}
  h\colon\Sigma\to E / G'
\end{displaymath}
of the bundle $E / G' = E\x_{G}(G/G')$, picking out a $G'$-orbit in each
fibre $E_x$. 

Let us now fix a maximal compact subgroup $H \subset G$. This choice,
together with the invariant inner product on $\lieg$, gives rise to a
\emph{Cartan decomposition}:
\begin{equation}
  \label{eq:5}
  \lie{g} = \lie{h} + \lie{m},
\end{equation}
where $\lieh$ is the Lie algebra of $H$ and $\liem$ is its orthogonal
complement. 

\begin{definition}
  A \emph{metric} in a principal $G$-bundle $E \to \Sigma$ is a
  reduction of structure group to $H \subset G$.
\end{definition}

In case $E_{\rho} = \widetilde{\Sigma}\x_{\rho}G$ is a flat bundle, we
have
\begin{displaymath}
  E / H =  \widetilde{\Sigma}\x_{\rho} (G/H),
\end{displaymath}
and hence a metric $h$ in $E$ corresponds to a
$\pi_1\Sigma$-equivariant map
\begin{displaymath}
  \tilde{h}\colon \widetilde{\Sigma} \to G/H.
\end{displaymath}

The energy of the metric $h$ is essentially the integral over $\Sigma$
of the norm squared of the derivative of $h$. In the following we make
precise this concept.  To start with we give $\Sigma$ a Riemannian
metric and note that $G/H$ is a Riemannian manifold. Hence we can
calculate the norm $\abs{D\tilde{h}(\tilde{x})}$ at any point
$\tilde{x}\in \wt{X}$. Furthermore, since the group $G$ acts on $G/H$
by isometries, the derivative of $\tilde{h}$ satisfies
\begin{displaymath}
  \abs{D\tilde{h}(\tilde{x})}
  = \abs{D\tilde{h}(\gamma\cdot\tilde{x})}
\end{displaymath}
for any $\gamma \in \pi_1\Sigma$. 
Alternatively, we may proceed as follows.  Let $T^vE\to \Sigma$ be the
vertical tangent bundle of $E$. The fact that $E$ is flat means that
there is a natural projection $p\colon TE\to T^vE$ and we can define
the vertical part of the derivative of $h$ as the composition
\begin{displaymath}
  Dh = p \circ dh \colon T\Sigma \to TE \to T^vE.
\end{displaymath}
Clearly we have $\abs{Dh(x)} = \abs{D\tilde{h}(\tilde{x})}$ for any
$\tilde{x}\in E_x$.

\begin{definition}
  Let $\Sigma$ be a closed oriented surface with a Riemannian metric
  and let $E \to \Sigma$ be a flat principal $G$-bundle. The
  \emph{energy} of a metric $h$ in $E$ is
  \begin{displaymath}
    \mathcal{E}(h) = \int_{\Sigma}\abs{Dh}^2\vol.
  \end{displaymath}
\end{definition}

\begin{remark}
  \label{rem:conformal-invariance}
  Recall that on a surface the integral of a one-form is conformally
  invariant. Hence it suffices to give $\Sigma$ a conformal structure
  in order to make the energy functional well defined. 
\end{remark}

\begin{definition}
  A metric $h$ in a flat $G$-bundle $E \to \Sigma$ is \emph{harmonic}
  if is a critical point of the energy functional.
\end{definition}

Next we want to reformulate this in terms of connections.  Let
$i\colon E_H \to E$ be the principal $H$-bundle obtained by the
reduction of structure group defined by the metric $h$, and denote the
flat connection on $E$ by $B \in \Omega^1(E,\lie{g})$. Using the
Cartan decomposition~(\ref{eq:5}) we can then write
\begin{equation}
  \label{eq:4}
  i^*B = A + \psi,
\end{equation}
where $A \in \Omega^1(E_H,\lie{h})$ defines a connection $E_H$ and
$\psi \in \Omega^1(E_H,\lie{m})$ is a tensorial $1$-form which
therefore descends to a section, abusively denoted by the same symbol,
\begin{displaymath}
  \psi \in \Omega^1(E_H(\lie{m})).
\end{displaymath}
Note that we have a canonical identification
\begin{displaymath}
  E_H(\lie{m}) \cong T^vE 
\end{displaymath}
and that under this identification we have 
\begin{displaymath}
  \psi = Dh,
\end{displaymath}
as is easily checked. To calculate the critical points of the energy
functional, take a deformation of the metric $h$ of the form
\begin{displaymath}
  h_t = \exp(t\cdot s)h \in \Omega^0(\Sigma,E/H)
\end{displaymath}
for $s \in \Omega^0(\Sigma,E_H(\lie{m}))$.
One then calculates
\begin{displaymath}
  \frac{d}{dt}(\mathcal{E}(h_t))_{|t=0}
  =\langle \psi,d_A s\rangle
\end{displaymath}
from which we deduce the following.
\begin{proposition}
  Let $h$ be a metric in a flat bundle $E\to \Sigma$ and let
  $(A,\psi)$ be defined by (\ref{eq:4}). Then  $h$ is harmonic if and
  only if
  \begin{displaymath}
    d_A^*\psi=0.
  \end{displaymath}
\end{proposition}


\subsection{The Corlette--Donaldson theorem}

The following result was proved independently by
Donaldson \cite{donaldson:1987} (for $G=\SL(2,\C)$) and Corlette
\cite{corlette:1988} (for more general groups and base manifolds
of dimension higher than two); see also Labourie \cite{labourie:1991}.  
The idea of the proof is to adapt the proof of  Eells--Sampson on the
existence of harmonic maps into negatively curved target manifolds to
the present ``twisted situation''.

\begin{theorem}
  \label{thm:donaldson-corlette}
  A flat bundle $E \to \Sigma$ corresponding to a representation
  $\rho\colon\pi_1\Sigma \to G$ admits a harmonic metric if and only
  if $\rho$ is reductive.
\end{theorem}

In terms of the pair $(A,\theta)$ given by (\ref{eq:4}), the flatness
condition on $B$ becomes
\begin{equation}
  \label{eq:7}
  \begin{aligned}
    F(A) + \tfrac{1}{2}[\theta,\theta] &=0, \\
    d_A\theta &=0,
  \end{aligned}
\end{equation}
as can be seen by considering the $\lieh$- and $\liem$-valued parts of
the equation $F(B)=0$ separately. This motivates the following
definition.

\begin{definition}
  Let $E_H\to X$ be a principal $H$-bundle on $X$, let $A$ be
  connection on $E_H$ and let
  $\theta\in\Omega^1(X,E_H(\liem))$. The triple $(E_H,A,\theta)$ is
  called a \emph{harmonic bundle} if the equations
    \begin{align}
      F(A) + \tfrac{1}{2}[\theta,\theta] &=0, \label{eq:8} \\
      d_A\theta &=0, \label{eq:12} \\
      d_A^*\theta &=0 \label{eq:13}
    \end{align}
  are satisfied.
\end{definition}


Next we want to obtain a statement at the level of moduli spaces
(analogous to (\ref{eq:21})). Fix a reduction $E_H \into E$ and
consider the gauge groups
\begin{displaymath}
  \begin{aligned}
    \mathcal{H} &= \Aut(E_H) = \Omega^0(\Sigma,E_H\x_{\Ad}H),\\
    \mathcal{G} &= \Aut(E) = \Omega^0(\Sigma,E\x_{\Ad}G).
  \end{aligned}
\end{displaymath}
Then Theorem~\ref{thm:donaldson-corlette} can equivalently be
formulated as saying that for any flat reductive connection $B$ in
$E$, there is a gauge transformation $g\in \mathcal{G}$ such that,
writing $g\cdot B = A + \psi$, the triple $(E_H,A,\psi)$ is a harmonic
bundle. 

Let $d \in\pi_1 H $ and fix $E_H$ with $c(E_H) = d $. The
\emph{moduli space of harmonic bundles of topological class $d$} is
\begin{displaymath}
  \Mhar_d(\Sigma,G) = \{(A,\theta)\suchthat \text{(\ref{eq:8})--(\ref{eq:13})
    hold}\} / \mathcal{H}. 
\end{displaymath}

Now Theorem~\ref{thm:donaldson-corlette} can be complemented by a
suitable uniqueness statement (analogous to
Proposition~\ref{prop:flat-rep}) which allows us to altogether obtain
a bijective correspondence
\begin{equation}
  \label{eq:DC-global}
  \MdR_d(\Sigma,G) \cong \Mhar_d(\Sigma,G).
\end{equation}


\section{Lecture 2: $G$-Higgs bundles and the Hitchin--Kobayashi correspondence}

In Lecture~1 we saw that any reductive surface group representation
gives rise to an essentially unique harmonic metric in the associated
flat bundle. In this lecture, we shall reinterpret this in holomorphic
terms, introducing  $G$-Higgs bundles. Moreover we shall explain the
Hitchin--Kobayashi correspondence for these.

Recall from Remark~\ref{rem:conformal-invariance} that we equipped the
surface $\Sigma$ with a conformal class of metrics. This is equivalent
to having defined a Riemann surface, which we shall henceforth denote
by $X = (\Sigma,J)$.

\subsection{Lie theoretic preliminaries}

Let $H^{\C}$ be the complexification of the maximal compact subgroup
$H \subseteq G$ and let $\liehc$ and $\liegc$ be the complexifications
of the Lie algebras $\lieh$ and $\lieg$, respectively. In particular,
$\liehc = \Lie(H^{\C})$. However we do not need to assume the
existence of a complexification of the Lie group $G$.

The Cartan decomposition (\ref{eq:5}) complexifies to
\begin{equation}
  \label{eq:6}
  \liegc = \liehc + \liemc;
\end{equation}
note that this is a direct sum of vector spaces but not of Lie
algebras. In fact, we have
\begin{displaymath}
  [\liehc,\liehc] \subseteq \liehc,\quad
  [\liehc,\liemc] \subseteq \liemc,\quad
  [\liemc,\liemc] \subseteq \liehc.
\end{displaymath}
Moreover, we have the $\C$-linear Cartan involution
\begin{displaymath}
  \theta\colon \liegc \to \liegc,
\end{displaymath}
whose $\pm 1$-eigenspace decomposition is (\ref{eq:6}), the 
real structure (i.e.\ $\C$-antilinear involution) corresponding to
$\lieg \subset \liegc$
\begin{displaymath}
  \sigma\colon \liegc \to \liegc
\end{displaymath}
and the compact real structure
\begin{displaymath}
  \tilde\tau=\theta\circ\sigma \colon \liegc \to \liegc.
\end{displaymath}
The $+1$-eigenspace of $\tilde\tau$ is a maximal compact subalgebra of
$\liegc$ whose intersection with  $\liehc$ is
$\lieh$.

We shall also need the \emph{isotropy representation} of $H^{\C}$ on
$\liemc$,
\begin{equation}
  \label{eq:isotropy-rep}
  \iota\colon H^{\C} \to \GL(\liemc),
\end{equation}
which is induced by the complexification of the adjoint
action of $H$ on $\lieg$

\subsection{The Hitchin equations}

We extend $\tilde\tau$ to
\begin{displaymath}
  \tau\colon \Omega^1(X,E(\liemc)) \to \Omega^1(X,E(\liemc))
\end{displaymath}
by combining it with conjugation on the form component. Locally
\begin{displaymath}
  \tau(\omega\otimes a) := \bar{\omega}\otimes\tau(a)
\end{displaymath}
for a complex $1$-form $\omega$ on $X$ and a section $a$ of
$E(\liemc)$. There is an isomorphism
\begin{equation}
  \label{eq:9}
  \begin{aligned}
    \Omega^1(E(\liem)) &\to \Omega^{1,0}(X,E(\liemc)), \\
    \theta &\mapsto \frac{\theta -iJ\theta}{2}
  \end{aligned}
\end{equation}
where $J$ is the complex structure on the tangent bundle of $X$. The
inverse given by
\begin{equation}
  \label{eq:10}
  \theta= \phi -\tau(\phi).
\end{equation}
This is entirely analogous to the way in which we can write the
connection $A \in \Omega^1(E_H,\lieh)$
\begin{equation}
  \label{eq:11}
  A = A^{1,0} + A^{0,1}
\end{equation}
with $A^{p,q}\in \Omega^{p,q}(E_H(\liemc))$.
\begin{remark}
  Note that
  \begin{displaymath}
    E_H(\liemc) = E_{H^{\C}}(\liemc),
  \end{displaymath}
  where $E_{H^{\C}} = E_H\x_H H^{\C}$ is the principal $H^{\C}$-bundle
  obtained by extension of structure group.  

  The bijective correspondence $A \leftrightarrow A^{0,1}$ gives us a
  bijective correspondence between connections $A$ on $E_H$ and
  holomorphic structures on $E_{H^{\C}}$ (the integrability condition
  is automatically satisfied because $\dim_{\C}X = 1$).

  Correspondingly, for any complex representation $V$ of $H^{\C}$, the
  vector bundle $E_{H^{\C}}(V)$ becomes a holomorphic bundle and the
  covariant derivative on sections of $E_{H^{\C}}(V)$ given by $A$
  decomposes as $d_A=\dbar_A+\partial_A$, where
  \begin{displaymath}
    \dbar_A\colon\Omega^0(X,E_{H^{\C}}(V)) \to \Omega^{0,1}(X,E_{H^{\C}}(V)).
  \end{displaymath}
  The holomorphic sections of $E_{H^{\C}}(V)$ are just the ones which
  are in the kernel of $\dbar_A$.
\end{remark}

With all this notation in place one sees, using the K\"ahler
identities, that the harmonic bundle equations
(\ref{eq:8})--(\ref{eq:13}) are equivalent to the \emph{Hitchin
  equations}
\begin{align}
  F(A) - [\phi,\tau(\phi)] &=0, \label{eq:14} \\
  \dbar_A\phi &=0. \label{eq:15}
\end{align}
Thus we have a canonical identification
\begin{equation}
  \label{eq:22}
  \Mhar_d(X,G) = \MHit_d(X,G),
\end{equation}
where we have introduced the moduli space
\begin{displaymath}
  \MHit_d(X,G)=\{(A,\phi)\suchthat \text{(\ref{eq:14})--(\ref{eq:15})
    hold}\} / \mathcal{H}
\end{displaymath}
of solutions to the Hitchin equations. This gauge theoretic point of
view allows one to give the moduli space $\MHit_d(X,G)$ a K\"ahler
structure. While the metric depends on the choice of conformal
structure on $\Sigma$, the K\"ahler form is independent of this choice,
and in fact coincides with Goldman's symplectic form \cite{goldman:1984}.

\subsection{$G$-Higgs bundles, stability and The Hitchin--Kobayashi cor\-respon\-dence}
\label{sec:g-higgs-bundles-HK}

The second Hitchin equation (\ref{eq:15}) says that $\Phi$ is
holomorphic with respect to the structure of holomorphic bundle. Write
$K = T^*X^{\C}$ for the holomorphic cotangent bundle, or
\emph{canonical bundle}, of $X$ and $H^0$ for holomorphic sections. We
have thus reached the conclusion that the harmonic bundle gives rise
to a holomorphic object, a so-called $G$-Higgs bundle, defined as
follows.

\begin{definition}
  A \emph{$G$-Higgs bundle} on $X$ is a pair $(E,\phi)$, where $E\to
  X$ is a holomorphic principal $H^{\C}$-bundle and $\phi \in
  H^0(X,E(\liemc)\otimes K)$.
\end{definition}

When $G$ is a complex group, we have that $H^{\C} = G$ and the Cartan
decomposition $\lieg^{\C} = \lieg + i \liegc$. Hence a $G$-Higgs
bundle is a pair $(E,\phi)$, where $E$ is a holomorphic principal
$G$-bundle and $\phi \in H^0(X,E(\lieg)\otimes K)$. Note that
$E(\liegc) = \Ad E$ is just the adjoint bundle of $E$. 

Another particular case is when $G = H$ is a compact group. Then we
have $\phi=0$, so a $G$-Higgs bundle is just a holomorphic principal
bundle and the Hitchin equations simply say that $F(A)=0$. 

In the following we give some examples of $G$-Higgs bundles for
specific groups. 

\begin{example}
  Let $G=\SU(n,\C)$. Then a $G$-Higgs bundle is just a holomorphic
  vector bundle $V \to X$ with trivial determinant.
\end{example}

\begin{example}
  Let $G=\SL(n,\C)$. Then a $G$-Higgs bundle is a pair $(V,\phi)$,
  where $V \to X$ is a holomorphic vector bundle with trivial
  determinant and $\phi\in H^0(X,\End_0(E)\otimes K)$ (where
  $\End_0(E)$ is the subspace of traceless endomorphisms).
\end{example}

\begin{example}
  Let $G = \SU(p,q)$. Then a $G$-Higgs bundle is a triple $(V,W,\phi)$,
  where $V$ and $W$ are holomorphic vector bundles on $X$ of rank $p$
  and $q$ respectively, satisfying
  $\det(V)\otimes\det(W)\cong\mathcal{O}$ and
  \begin{displaymath}
    \phi = (\beta,\gamma) \in H^0(X,\Hom(W,V)\otimes K)
      \oplus H^0(X,\Hom(V,W)\otimes K).
  \end{displaymath}
\end{example}

\begin{example}
  Let $G = \Sp(2n,\R)$. Then a $G$-Higgs bundle is a pair $(V,\phi)$,
  where $V$ is a holomorphic vector bundle on $X$ of rank $n$ and
  \begin{displaymath}
    \phi = (\beta,\gamma) \in H^0(X,S^2V \otimes K)
      \oplus H^0(X,S^2V^*\otimes K).
  \end{displaymath}
\end{example}

In case $G=\SU(n)$, we are thus in the presence of a complex vector
bundle with a flat unitary connection. Such a bundle turns out to be
\emph{polystable}. The Narasimhan--Seshadri Theorem
\cite{narasimhan-seshadri:1965}, conversely, says that if a holomorphic
vector bundle is polystable then it admits a metric such that the
unique unitary connection compatible with the holomorphic structure is
(projectively) flat. There is an analogous statement for other
compact $G$, due to Ramanathan \cite{ramanathan:1975}).

These results generalize to $G$-Higgs bundles. The appropriate
stability condition is a bit involved to state in general. However, in
the case of Higgs vector bundles it is simply the following. Recall
that the \emph{slope} of a vector bundle $E \to X$ is $\mu(E) =
\deg(E)/\rk(E)$. Also, we say that a subbundle $F \subset E$ is
\emph{$\phi$-invariant} if $\phi(F) \subset F \otimes K$.

\begin{definition}
  A Higgs vector bundle $(E,\phi)$ is \emph{semistable} if
  \begin{displaymath}
    \mu(F) \leq \mu(E)
  \end{displaymath}
  for any subbundle $\phi$-invariant subbundle $F \subset E$ and it is
  \emph{stable} if, moreover, strict inequality holds whenever $F$ is
  proper and non-zero. A Higgs vector bundle $(E,\phi)$ is
  \emph{polystable} if it is isomorphic to a direct sum of stable
  Higgs bundles, all of the same slope.
\end{definition}

The stability conditions for $G$-Higgs bundles can be obtained as a
special case of a general stability conditions for \emph{pairs} and we
refer the reader to \cite{garcia-gothen-mundet:2009b} for the detailed
formulation. It is worth noting that
poly- and semistability of a
$G$-Higgs bundle $(E,\phi)$ are equivalent to poly- and semistability
of the Higgs vector bundle $(E(\liegc),\ad(\phi))$.

The general Hitchin--Kobayashi correspondence for principal pairs
\cite{mundet:2000,bradlow-garcia-prada-mundet:2003} now has as a
consequence the following Hitchin--Kobayashi correspondence for
$G$-Higgs bundles. (see \cite{garcia-gothen-mundet:2009b} for the
full extension to polystable pairs, as well as a detailed analysis of
the case of $G$-Higgs bundles.)

\begin{theorem}
  \label{thm:G-higgs-HK}
  Assume that $G$ is semisimple. A $G$-Higgs bundle $(E,\phi)$ is
  polystable if and only if it admits a reduction of structure group
  to the maximal compact $H \subset H^{\C}$, unique up to isomorphism
  of $H$-bundles, such that the following holds: denoting by $A$ the
  unique $H$-connection compatible with the reduction and $\dbar_A$
  the $\dbar$-operator induced form the holomorphic structure, the
  pair $(A,\phi)$ satisfies the Hitchin equations (\ref{eq:14}) and
  (\ref{eq:15}).
\end{theorem}

In the context of Higgs bundles, the Hitchin--Kobayashi correspondence
goes back to the work of Hitchin \cite{hitchin:1987a} and Simpson
\cite{simpson:1988}.

\begin{remark}
  The assumption that $G$ be semisimple is not essential. Indeed, as
  can be expected from the situation for usual $G$-bundles, for
  reductive $G$ an analogous statement holds, if one adds a suitable
  central term to the right hand side of the first of the Hitchin
  equations. For the correspondence with representations, one must
  then consider homomorphisms of a central extension of the
  fundamental group. We refer to \cite{garcia-gothen-mundet:2009b} for
  more details on this.
\end{remark}

Just as for vector bundles, stability of $G$-Higgs bundles has a dual
importance. Namely, apart from its role in the Hitchin--Kobayashi
correspondence, it is also the appropriate notion for constructing moduli
spaces using GIT. The constructions Schmitt (see the book
\cite{schmitt:2008}) are in fact sufficiently general to also cover
many cases of $G$-Higgs bundles. Thus we have yet another moduli space
at our disposal, namely the moduli space
\begin{displaymath}
  \MDol_d(X,G)
\end{displaymath}
of semistable $G$-Higgs bundles of topological class $d \in \pi_1 H$.
Alternatively, this moduli space can be constructed using a Kuranishi
slice method. From this point of view, we fix a principal $H^{\C}$-bundle
$E\to \Sigma$ and consider the complex configuration space
\begin{displaymath}
  \mathcal{C}^{\C} = \{(A^{0,1},\phi) \suchthat \dbar_A\phi=0\}.
\end{displaymath}
The complex gauge group $\mathcal{H}^{\C}$ acts naturally on this
space and on the subspace $\mathcal{C}^{\C}_{\text{polystable}}$ of
pairs $(A^{0,1},\phi)$ which define the structure of a polystable $G$-Higgs
bundle on $E$. The moduli space is then
\begin{displaymath}
  \MDol_d(X,G) = \mathcal{C}^{\C}_{\text{polystable}} / \mathcal{H}^{\C}.
\end{displaymath}
Either way, Theorem~\ref{thm:G-higgs-HK} implies that we have an
identification 
\begin{equation}
  \label{eq:23}
  \MDol_d(X,G) \cong \MHit_d(X,G).
\end{equation}
Putting together this with the previous identifications (\ref{eq:21}),
(\ref{eq:DC-global}) and (\ref{eq:22}) we finally obtain the
\emph{non-abelian Hodge Theorem}.

\begin{theorem}
  \label{thm:non-ab-Hodge}
  Let $X$ be a closed Riemann surface of genus $g$. Then there is a
  homeomorphism
  \begin{displaymath}
    \MB_d(X,G) \cong \MDol_d(X,G).
  \end{displaymath}
\end{theorem}

\begin{remark}
  The fact that the identification of Theorem~\ref{thm:non-ab-Hodge}
  is a homeomorphism is not too hard to see, but more is true: outside
  of the singular loci, the identification is in fact an analytic
  isomorphism. On the other hand, it is definitely not algebraic. In
  this respect, it is instructive to consider the example $G=\C^*$.
\end{remark}

\subsection{The Hitchin map}

We end this Lecture by recalling the definition of the Hitchin map
which plays a central role in the theory of Higgs bundles.

Take a basis $\{p_1,\dots,p_r\}$ for the invariant polynomials on the
Lie algebra $\liegc$ and let $d_i = \deg(p_i)$. Given a $G$-Higgs
bundle $(E,\phi)$, evaluating $p_i$ on $\phi$ gives a section
$p_i(\phi) \in H^0(X,K^{d_i})$. The \emph{Hitchin map} is defined to be
\begin{equation}
  \label{def-hitchin-map}
  \begin{aligned}
    h\colon \MDol_d &\to B, \\
    (E,\phi) & \mapsto (p_1(\phi),\dots,p_r(\phi)),
  \end{aligned}
\end{equation}
where the \emph{Hitchin base} is
\begin{displaymath}
  B := \bigoplus H^0(X,K^{d_i}).
\end{displaymath}
The Hitchin map is proper and, for $G$ complex, defines an
algebraically completely integrable system known as the Hitchin system
(see Hitchin \cite{hitchin:1987b}.)

\subsection{The moduli space of $\SU(p,q)$-Higgs bundles}
\label{sec:supq}

We end this section by illustrating how the Higgs bundle point of view
allows for easy proofs of strong results by proving the Milnor--Wood
inequality for $\SU(p,q)$-Higgs bundles, and discussing a closely
related rigidity result.

Recall that an $\SU(p,q)$-Higgs bundle is a quadruple
$(V,W,\beta,\gamma)$, where $V$ and $W$ are vector bundles on $X$ of
rank $p$ and $q $ respectively, satisfying
$\det(V)\otimes\det(W)\cong\mathcal{O}$, and where $\beta\in
H^0(X,\Hom(W,V)\otimes) K $ and $\gamma\in H^0(X,\Hom(V,W)\otimes
K)$. The topological classification of such bundles is given by
$\deg(V) = -\deg(W) \in \Z$. Denote by $\mathcal{M}_d$ the moduli space of
$\SU(p,q)$-Higgs bundles with $\deg(V)=d$.

In the case $p=q=1$, we have $\SU(1,1) = \SL(2,\R)$ and the degree $d$
is just the Euler class of the corresponding flat
$\SL(2,\R)$-bundle. In 1957 Milnor \cite{milnor:1957} proved that it
satisfies the bound
\begin{displaymath}
  \abs{d}\leq g-1. 
\end{displaymath}
Much more generally, whenever $G$ is non-compact of Hermitian type, one
can define an integer invariant, the \emph{Toledo invariant}, of
representations $\rho\colon\pi_1 \Sigma \to G$ and there is a bound on
the Toledo invariant, usually known as a \emph{Milnor--Wood
  inequality}. In various degrees of generality this is due to, among
others, Domic--Toledo \cite{domic-toledo:1987}, Dupont
\cite{dupont:1978}, Toledo \cite{toledo:1979,toledo:1989} and Turaev
\cite{turaev:1984}. In the case of $G=\SU(p,q)$, the Milnor--Wood
inequality is
\begin{equation}
  \label{eq:supq-MW}
  \abs{d} \leq \min\{p,q\}(g-1).
\end{equation}
A proof of this Milnor--Wood inequality using Higgs bundles is very
easy to give. For this it is convenient (though not essential) to pass
through usual Higgs vector bundles: since $\SU(p,q)$ is a subgroup
of $\SL(p+q,\C)$, to any $\SU(p,q)$-Higgs bundle we can associate an
$\SL(p+q,\C)$-Higgs bundle
\begin{displaymath}
  (E,\Phi) = (V\oplus W,
  \begin{pmatrix}
    0 & \beta \\
    \gamma & 0
  \end{pmatrix}
).
\end{displaymath}
Now, if $(V,W,\beta,\gamma)$ is polystable, then so is $(E,\Phi)$ (this
follows immediately from the fact that a solution to the
$\SU(p,q)$-Hitchin equations on $(V,W,\beta,\gamma)$ induces a
solution on $(E,\Phi)$). Let $N \subset V$ be the kernel of
$\gamma\colon V \to W\otimes K$, viewed as a subbundle, and let $I
\otimes K
\subset W \otimes K$ be the subbundle obtained by saturating the image
subsheaf. Thus, $\gamma$ induces a bundle  map of maximal rank
$\bar{\gamma}\colon V/N
\to I \otimes K$, from which we deduce that
\begin{equation}
  \label{eq:17}
  \deg(N) - \deg(V) + \deg(I) + (2g-2)\rk(\gamma) \geq 0 
\end{equation}
with equality if and only if $\bar{\gamma}$ is an isomorphism.  
Moreover, the subbundles $N \subset E$ and $V \oplus I \subset E$ are
$\Phi$-invariant, so polystability of $(E,\Phi)$ implies  that 
\begin{align}
  \deg(N) &\leq 0, \label{eq:18} \\
  \deg(V) + \deg(I) &\leq 0. \label{eq:19}
\end{align}
Putting together equations (\ref{eq:17})--(\ref{eq:19}) we obtain
\begin{equation}
  \label{eq:20}
  \deg(V) \leq \rk(\gamma)(g-1)
\end{equation}
from which the Milnor--Wood inequality (\ref{eq:supq-MW}) is
immediate for $d \geq 0$. When $d\leq 0$ a similar argument involving
$\beta$ instead of $\gamma$ gives the result. 

But our arguments in fact give more information in the case
when equality holds in (\ref{eq:supq-MW}). Assume for definiteness
that $d \geq 0$ and that $p \leq q$. Then, if equality holds in
(\ref{eq:supq-MW}), we conclude immediately that $\rk(\gamma) = p$ and
that $\gamma\colon V \to I \otimes K$ is an isomorphism. Hence, by
polystability of $(E,\Phi)$, there is a decomposition $W = I \oplus Q$
and $\beta_{|Q}=0$. In other words, the $\SU(p,q)$-Higgs bundle
$(V,W,\beta,\gamma)$ decomposes into the $\U(p,p)$-Higgs bundle
$(V,I,\beta,\gamma)$ and the $\U(q-p)$-Higgs bundle $Q$.

From the point of view of representations of the fundamental group
this can be viewed as a rigidity result which was first proved by
Toledo \cite{toledo:1989} for $p=1$ and by Hern\'andez
\cite{hernandez:1991} for $p=2$. A more general result valid in the
context of arbitrary groups of Hermitian type has been proved
Burger--Iozzi--Wienhard
\cite{burger-iozzi-wienhard:2010,burger-iozzi-wienhard:2003}. From the
point of view of Higgs bundles, the results for $\U(p,q)$ appeared in
\cite{bradlow-garcia-prada-gothen:2003} and a survey of the situation
for other classical groups can be found in
\cite{bradlow-garcia-prada-gothen:hss-higgs}, while the PhD thesis of Rubio
\cite{rubio:2012} treats the question for general groups using a
general Lie theoretic approach.

\section{Lecture 3: Morse--Bott theory of the moduli space of $G$-Higgs bundles}

In this final lecture we consider the $\C^*$-action on the moduli
space of $G$-Higgs bundles and explain how to use it to study its
topology. We shall consider the Dolbeault moduli space and
occasionally use the identification with the gauge theory moduli space
of solutions to Hitchin's equation. For simplicity we shall denote it
simply by $\mathcal{M}_d$. Again, though not strictly necessary, we
shall assume that $G$ is semisimple. To get started we need to review
some of the deformation theory of $G$-Higgs bundles.

\subsection{Simple and infinitesimally simple $G$-Higgs bundles}

Let $E$ be a principal $H^{\C}$-bundle on $X$. An \emph{automorphism}
of $E$ is an equivariant holomorphic bundle map $g\colon E \to E$
which admits a holomorphic inverse. We denote the group of
automorphisms of $E$ by $\Aut(E)$.
Equivalently, we may define $\Aut(E)$ to be the space of holomorphic
sections of the bundle of automorphisms $E\times_{\Ad}G\to X$.
Let $(E,\phi)$ be a $G$-Higgs bundle. We denote by $\Aut(E)$ the group
of automorphisms of $(E,\phi)$:
\begin{displaymath}
  \Aut(E,\phi) = \{g \in \Aut(E) \suchthat \Ad(g)(\phi)=\phi \}.
\end{displaymath}
We also introduce the \emph{infinitesimal automorphism space} (which,
at least formally, is the Lie algebra of the automorphism group), defining
\begin{displaymath}
  \aut(E,\phi) = \{Y \in H^0(X,E(\liehc)) \suchthat [Y,\phi] = 0\}.
\end{displaymath}

A $G$-Higgs bundle $(E,\phi)$ is \emph{simple} if its automorphism
group is smallest possible, i.e.,
\begin{displaymath}
  \Aut(E,\phi) = Z(H^{\C})\cap\ker(\iota).
\end{displaymath}
Also, we say that $(E,\phi)$ is \emph{infinitesimally simple} if
\begin{displaymath}
  \aut(E,\phi) = Z(\liehc)\cap\ker(d\iota).
\end{displaymath}
Note that for Higgs vector bundles, these two notions are
equivalent. This is, however, not true in general, as
Example~\ref{ex:stable-not-simple} below shows.

The following result is the $G$-Higgs bundle version of the well known
fact that a stable vector bundle only has scalar automorphisms.

\begin{proposition}
  Let $(E,\phi)$ be a stable $G$-Higgs bundle. Then it is
  infinitesimally simple.
\end{proposition}


\begin{example}
  \label{ex:stable-not-simple}
  Let $M_1$ and $M_2$ be line bundles on $X$ with $M_i^2=K$ and $M_1
  \neq M_2$. Define $V = M_1 \oplus M_2$ and let $\beta=0\in H^0(X,
  S^2V\otimes K)$ and 
  $\gamma=\left(
  \begin{smallmatrix}
    1 & 0 \\
    0 & 1
  \end{smallmatrix}\right)\in H^0(X,S^2V^*\otimes K)$.
  Then it is easy to see that the $\Sp(2,\R)$-Higgs bundle $(V,\beta,\gamma)$ is
  stable and hence infinitesimally simple. However, it is not simple
  since it has the automorphism $\left(
  \begin{smallmatrix}
    -1 & 0 \\
    0 & 1
  \end{smallmatrix}\right)$.
\end{example}

\subsection{Deformation theory of $G$-Higgs bundles}
\label{sec:deformation-theory}

Next we outline the deformation theory of $G$-Higgs bundles.  A useful
reference for the following material is Biswas--Ramanan
\cite{biswas-ramanan:1994}.

\begin{definition}
  The \emph{deformation complex} of a $G$-Higgs bundle $(E,\phi)$ is
  the complex of sheaves
  \begin{displaymath}
    C^\bullet(E,\phi)\colon
    E(\liehc)\xra{[-,\phi]}E(\liemc)\otimes K.
  \end{displaymath}
\end{definition}

The deformation theory of a $G$-Higgs bundle $(E,\phi)$ is governed by
the hypercohomology groups of the deformation complex. Thus, we have
the following standard results.

\begin{proposition}
  Let $(E,\phi)$ be a $G$-Higgs bundle.
  \begin{enumerate}
  \item There is a canonical identification between the space of
    infinitesimal deformations of $(E,\phi)$ and the hypercohomology
    group
    \begin{displaymath}
      \HH^1(C^\bullet(E,\phi))
    \end{displaymath}

  \item There is a long exact sequence 
\begin{displaymath}
  \begin{split}
  0 &\to \HH^0(C^{\bullet}(E,\varphi)) \to H^0 (E(\lieh^\CC))
  \xrightarrow{[-,\phi]} H^0(E(\liem^\CC)\otimes K) \\
  &\to \HH^1(C^{\bullet}(E,\varphi))
  \to H^1 (E(\lieh^\CC)) \xrightarrow{[-,\phi]}
  H^1(E(\liem^\CC)\otimes K) \\
  &\to
  \HH^2(C^{\bullet}(E,\varphi)) \to 0.
  \end{split}
\end{displaymath}
\end{enumerate}
\end{proposition}

Note that the long exact sequence in (2) of the previous Proposition
immediately implies that there is a canonical identification
\begin{displaymath}
  \aut(E,\Phi) = \HH^0(C^\bullet(E,\phi)).
\end{displaymath}

One way of proving the following result is to consider the Kuranishi
slice method for constructing the moduli space mentioned in
Section~\ref{sec:g-higgs-bundles-HK}.

\begin{proposition}
  \label{prop:smoothness}
  Assume that $(E,\phi)$ is a stable and simple $G$-Higgs bundle and
  that the vanishing $\HH^2(C^\bullet(E,\phi)) = 0$ holds. Then
  $(E,\phi)$ represents a smooth point of the moduli space $\mathcal{M}_d$.
\end{proposition}

\begin{remark}
  If $G$ is a reductive group which is not necessarily semisimple, one
  should consider the \emph{reduced deformation complex}, obtained by
  dividing out by $Z(\liegc)$; equivalently, this is the deformation
  complex of the $PG$-Higgs bundle obtained from the $G$-Higgs bundle.
\end{remark}

\subsection{The $\C^*$-action and topology of moduli spaces}

In order to avoid the problems arising from the presence of
singularities, throughout this section we shall make the assumption
that we are in a situation where the moduli space $\mathcal{M}_d$ is
smooth. 

It is a very important feature of the moduli space of Higgs bundles
that it admits an action of the multiplicative group of
non-zero complex numbers:
\begin{equation}
  \label{eq:C-star-action}
  \begin{aligned}
    \C^*\x \mathcal{M}_d &\to \mathcal{M}_d \\
    (z, (E,\phi)) &\mapsto (E,z\phi).
  \end{aligned}
\end{equation}

There are two distinct ways of using this action to obtain topological
information about the moduli space, as we shall now explain. However,
a theorem of Kirwan ensures that they give essentially equivalent
information. 

We start by a Morse theoretic point of view. For this we use the
identification between the Dolbeault moduli space and the moduli space
of solutions to the Hitchin equations \eqref{eq:11}--\eqref{eq:14} given by
Theorem~\ref{thm:G-higgs-HK}. Observe that the subgroup $S^1
\subset \C^*$ acts on the moduli space of solutions to the
Hitchin equations. With respect to the
(symplectic) K\"ahler form on the moduli space this action is
Hamiltonian and it has a moment map (up to some fixed scaling) given by
\begin{align*}
  f\colon\mathcal{M}_d &\to \R \\
  (A,\phi) &\mapsto \norm{\phi}^2:=\int_X\abs{\phi}^2\vol.
\end{align*}
Hitchin \cite{hitchin:1987a} showed, using Uhlenbeck's weak
compactness theorem, that $f$ is a proper map.  Moreover, it follows
from a theorem of Frankel \cite{frankel:1959} that $f$ is a perfect
Bott--Morse function. That $f\colon M \to \R$ is \emph{Bott--Morse}
means that its critical points form smooth (connected) submanifolds
$N_\lambda \subseteq M$ such that the Hessian of $f$ is non-degenerate
along the normal bundle to $N_\lambda$ in $M$. That $f$ is
\emph{perfect} means that the Poincar\'e polynomial
\begin{displaymath}
  P_t(M):= \sum t^i \dim H^i(M,\QQ)
\end{displaymath}
can be determined as
\begin{equation}
  \label{eq:poincare-perfect}
  P_t(M) = \sum_{\lambda} t^{\mathrm{Index}(N_\lambda)}P_t(N_\lambda);
\end{equation}
here $\mathrm{Index}(N_\lambda)$ is the \emph{index} of the critical
submanifold $N_\lambda$, i.e., the real rank of the subbundle of the normal
bundle on which the Hessian of $f$ is negative definite.

The condition for $f$ to be a moment map for the Hamiltonian
$S^1$-action on $\mathcal{M}_d$ is
\begin{displaymath}
  \grad f = i\xi,
\end{displaymath}
where $\xi$ is the vector field generating the $S^1$-action. In
particular, the critical submanifolds of $f$ are just the components
of the fixed locus of the $S^1$-action. Moreover, if we denote by
$N^+_\lambda$ the stable manifold of $N_\lambda$, we obtain a
\emph{Morse stratification}
\begin{displaymath}
  \mathcal{M}_d = \bigcup_{\lambda} N^+_\lambda. 
\end{displaymath}
Note that the fact that $f$ is proper and bounded below guarantees
that every point in $\mathcal{M}_d$ belongs to one of the $N^+_\lambda$

The more algebraic point of view comes about by looking at the full
$\C^*$-action on $\MDol_d$. It is a general result of Bia{\l}ynicki-Birula
\cite{bialynicki-birula:1973} that there is an algebraic stratification defined
as follows: let $\{\tilde{N}_\lambda\}$ be the
components of the fixed locus and define
\begin{displaymath}
  \tilde{N}_\lambda^+ = \{m\in \mathcal{M}_d \suchthat \lim_{z\to 0}
  z\cdot m \in \tilde{N}_\lambda\}.
\end{displaymath}
Then the \emph{Bia{\l}ynicki-Birula stratification} is
\begin{displaymath}
  \mathcal{M}_d = \bigcup_\lambda \tilde{N}^+_\lambda.
\end{displaymath}
It is perhaps not immediately clear that every point in
$\mathcal{M}_d$ lie in one of the $\tilde{N}^+_\lambda$. It follows,
however, from the properness and equivariance (with respect to the
suitable weighted $\C^*$-action on the Hitchin base $B$) of the
Hitchin map (\ref{def-hitchin-map}).

The whole picture fits into the general setup of $\C^*$-actions on
K\"ahler manifolds arising from hamiltonian circle actions. In
particular, it follows from the results of Kirwan \cite{kirwan:1984}
that the Morse and Bia{\l}ynicki-Birula stratifications coincide.

From either point of view, one can now obtain topological information
on the moduli space, as pioneered by Hitchin \cite{hitchin:1987a} in
his calculation of the Poincar\'e polynomial of the moduli space of
the moduli space of rank 2 Higgs bundles. In general, the success of
this approach depends crucially on having a good understanding of the
topology of the fixed loci $N_\lambda$. We remark that the role played
by the underlying geometric decomposition of the moduli space is
perhaps best brought out by studying in the first place the class of
the spaces under study in the $K$-theory of varieties, and then
obtaining from this information such as Hodge and Poincar\'e
polynomials. For examples of this point of view we refer to
Chuang--Diaconescu--Pan \cite{Chuang:2010ii} or
Garc\'\i{}a-Prada--Heinloth--Schmitt
\cite{garcia-heinloth-schmitt:2011}.

\subsection{Calculation of Morse indices}
\label{sec:morse-indices}

Let us consider the fixed points of the circle action on
$\mathcal{M}_d$. For simplicity we start out with an ordinary Higgs vector
bundle $(E,\phi)$, where $E$ is a vector bundle and $\phi \in
H^0(X,\End(E)\otimes K)$.

The following is easily proved (see Hitchin \cite{hitchin:1987a} or
Simpson \cite{simpson:1992}).

\begin{proposition}
  The Higgs bundle $(E,\phi)$ is a fixed point of the circle action on
  $\MDol_d$ if and only if it is a \emph{Hodge bundle}, i.e., there is
  a decomposition
  \begin{displaymath}
    E = E_0 \oplus \dots E_p
  \end{displaymath}
  and, with respect to this decomposition, $\phi$ has weight one, by
  which we mean that $\phi(E_k) \subseteq E_{k+1}\otimes K$.
\end{proposition}

The basic idea is that the weight $k$ subbundle $E_k \subset E$ is the
$ik$-eigenbundle of the infinitesimal automorphism $\psi =
\lim_{\theta\to 0}g(\theta)$ counteracting
the circle action, where
\begin{displaymath}
  (E,e^{i\theta}\phi) = g(\theta)\cdot(E,\phi).
\end{displaymath}

For $G$-Higgs bundles in general, the simplest procedure is to work
out the shape of the Hodge bundles (fixed under the circle action) in
each individual case. Note that if the $G$-Higgs bundle $(E,\phi)$ is
fixed, then so is the the adjoint Higgs vector bundle
$(E(\liegc),\ad(\phi))$ and therefore it is a Hodge bundle. Moreover,
since the infinitesimal automorphism $\psi$ lies in $E(\lieh)$, the
decomposition of $E(\liegc)$ in eigenbundles is compatible with the
decomposition $E(\liegc) = E(\liehc) \oplus E(\liemc)$. It follows
that there are decompositions
\begin{displaymath}
  \begin{aligned}
    E(\liehc) &= \bigoplus E(\liehc)_k, \\
    E(\liemc) &= \bigoplus E(\liemc)_k,
  \end{aligned}
\end{displaymath}
and that with respect to these we have
\begin{displaymath}
  \begin{aligned}
    \ad(\phi)\colon E(\liehc)_k \to E(\liemc)_{k+1}\otimes K, \\
    \ad(\phi)\colon E(\liemc)_k \to E(\liehc)_{k+1}\otimes K.
  \end{aligned}
\end{displaymath}
In particular, the deformation complex of $(E,\phi)$ decomposes as
\begin{displaymath}
  C^\bullet(E,\phi) = \bigoplus C_k^\bullet(E,\phi), 
\end{displaymath}
where the weight $k$ piece of the deformation complex is given by
\begin{equation}
  \label{eq:16}
  C_k^\bullet(E,\phi)\colon E(\liehc)_k \xra{[-,\phi]}
  E(\liemc)_{k+1}\otimes K.
\end{equation}

An easy calculation (see for example \cite{garcia-gothen-munoz:2004}
for the case of ordinary parabolic Higgs bundles which is essentially
the same as the present one) now shows the following.

\begin{proposition}
  Let $(E,\phi)$ be a stable $G$-Higgs bundle which is fixed under the circle
  action and represents a smooth point of the moduli
  space. With the notation introduced above, we have
  \begin{align}
    \dim N^+_\lambda &= \dim \HH^1(C^\bullet_{\leq 0}), \\
    \dim N^-_\lambda &= \dim \HH^1(C^\bullet_{\geq 0}), \\
    \dim N_\lambda &= \dim \HH^1(C^\bullet_{0}).
  \end{align}
  Hence the Morse index of the critical submanifold
  $N_\lambda$ is
  \begin{displaymath}
    \Index(N_\lambda) = 2 \dim \HH^1(C^\bullet_{> 0}).
  \end{displaymath}
\end{proposition}

Bott--Morse theory shows that the number of 
connected components of the moduli space equals that of the subspace
of local minima of $f$. Thus, for the determination of this most basic
of topological invariants it is important to have a convenient
criterion for the Morse index to be zero. This is provided by the
following result
(\cite[Proposition~4.14]{bradlow-garcia-prada-gothen:2003}; see
\cite[Lemma~3.11]{bradlow-garcia-prada-gothen:homotopy} for a
corrected proof).

\begin{proposition}
  \label{prop:morse-zero-criterion}
  Let $(E,\phi)$ represent a critical point of $f$. Then $(E,\phi)$
  represents a local minimum of $f$ if and only if the map
  \begin{displaymath}
    [-,\phi]\colon E(\liehc) \to E(\liemc) \otimes K
  \end{displaymath}
  is an isomorphism for all $k>0$.
\end{proposition}

\subsection{The moduli space of $\Sp(2n,\R)$-Higgs bundles}

We end by illustrating how the ideas explained in this section work,
by considering the case $G =\Sp(2n,\R)$.

Recall that an $\Sp(2n,\R)$-Higgs bundle is a triple
$(V,\beta,\gamma)$, where $V\to X$ is a rank $n$ vector bundle,
$\beta\in H^0(X,S^2V)\otimes K $ and $\gamma\in H^0(X,S^2V^*)$. The
topological classification of such bundles is given by $\deg(V) \in
\Z$. Denote by $\mathcal{M}_d$ the moduli space of $\Sp(2n,\R)$-Higgs
bundles with $\deg(V)=d$. 

In the following we outline the application of the Morse theoretic
point of view for determining the number of connected components of
$\mathcal{M}_d$. We should point out that $\mathcal{M}_d$ is not a
smooth variety, so that care must be taken in dealing with
singularities in applying the theory. We shall ignore this issue for
reasons of space, and in order to bring out more clearly the main
ideas. We refer to \cite{garcia-gothen-mundet:2009a} for full details. 

Note that a $\Sp(2n,\R)$-Higgs bundle is in particular an
$\SU(n,n)$-Higgs bundle. Hence we have from (\ref{eq:supq-MW}) that
the Milnor--Wood inequality
\begin{equation}
  \label{eq:sp2nr-MW}
  \abs{d} \leq n(g-1)  
\end{equation}
holds. Say that a $\Sp(2n,\R)$-Higgs bundle is \emph{maximal} if
equality holds. 

Note that taking $V$ to its dual and interchanging $\beta$ and
$\gamma$ defines an isomorphism $\mathcal{M}_d \cong
\mathcal{M}_{-d}$. Hence we shall assume without loss of generality
that $d \geq 0$ for the remainder of this section.

Denote by $N_0 \subset \mathcal{M}_d$ the subspace of local minima of
$f$. In the non-maximal case,
Proposition~\ref{prop:morse-zero-criterion} leads to the following
result.
\begin{proposition}
  \label{prop:sp2nr-minima-non-max}
  Assume that $0 < d < n(g-1)$. Then the subspace of local minima
  $N_0\subset \mathcal{M}_d$ consists of all $(V,\beta,\gamma)$ with
  $\beta=0$. If $d=0$, the subspace of local minima $N_0\subset
  \mathcal{M}_0$ consists of all $(V,\beta,\gamma)$ with $\beta=0$ and
  $\gamma=0$.
\end{proposition}

Thus for $d=0$, the subspace $N_0$ can be identified with the moduli
space of polystable vector bundles of degree zero. Since this moduli
space is known to be connected, we conclude that $\mathcal{M}_0$ is
also connected.

For $0<d<n(g-1)$, the moduli space $\mathcal{M}_d$ is known to be
connected only for $n=1$ (by the results of Goldman
\cite{goldman:1988}, reproved by Hitchin \cite{hitchin:1987a} using
Higgs bundles) and for $n=2$ by Garc\'\i{}a-Prada--Mundet
\cite{garcia-prada-mundet:2004} (see also
\cite{gothen-oliveira:2012}). However, for $n \geq 3$, the
connectedness of $N_0$ --- and hence $\mathcal{M}_d$ --- appears to be
difficult to establish.  

On the other hand, when $d=n(g-1)$ is maximal, the complete answer is
known from the work of Goldman and Hitchin cited above when
$n=1$, from \cite{gothen:2001} when $n=2$, and from
\cite{garcia-gothen-mundet:2009a} when $n \geq 3$. It is as follows.

\begin{theorem}
  Let $\mathcal{M}_{\mathrm{\max}}$ be the moduli space of
  $\Sp(2n,\R)$-Higgs bundles $(V,\beta,\gamma)$ with
  $\deg(V)=n(g-1)$. Then
  \begin{enumerate}
  \item $\#\pi_0\mathcal{M}_{\mathrm{\max}}=2^{2g}$ for $n=1$,
  \item $\#\pi_0\mathcal{M}_{\mathrm{\max}}=3\cdot2^{2g}+2g-4$ for
    $n=2$, and
  \item $\#\pi_0\mathcal{M}_{\mathrm{\max}}=3\cdot2^{2g}$ for
    $n\geq 3$.
  \end{enumerate}
\end{theorem}

We end by briefly explaining how this result comes about. Hitchin
\cite{hitchin:1992} showed that whenever $G$ is a split real form, the
moduli space of $G$-Higgs bundles has a distinguished component, now
known as the \emph{Hitchin component}, which can be concisely
described in terms of representations of the fundamental group: it 
consists of $G$-Higgs bundles corresponding to representations which
factor through a Fuchsian representation of the fundamental group in
$\SL(2,\R)$, where $\SL(2,\R) \into G$ is embedded as a so-called
principal three-dimensional subgroup (when $G=\Sp(2n,\R)$ this is just
the irreducible representation of $\SL(2n,\R)$ on $\R^{2n}$). For
every $n$, the moduli space $\mathcal{M}_{\max}$ has
$2^{2g}$ Hitchin components which, however, become identified if the
pass to the projective group $\mathrm{PSp}(2n,\R)$.

To explain the appearance of the remaining components, recall the
argument used to prove the Milnor--Wood inequality (\ref{eq:supq-MW})
in Section~\ref{sec:supq}. This shows that for a maximal
$\Sp(2n,\R)$-Higgs bundle $(V,\beta,\gamma)$ (with $d\geq 0$), we have
an isomorphism
\begin{displaymath}
  \gamma\colon V \to V^*\otimes K.
\end{displaymath}
Hence, since $\gamma$ is symmetric, $V$ admits a $K$-valued everywhere
non-degenerate quadratic form. Defining $W=V \otimes K^{-n/2}$ and
$Q=\gamma\otimes 1_{K^{-n/2}}$ we obtain an $\mathrm{O}(n,\C)$-bundle
$(W,Q)$, meaning that we obtain new topological invariants defined by
the Stiefel--Whitney classes $w_1$ and $w_2$ of $(W,Q)$. These then
give rise to new subspaces $\mathcal{M}_{w_1,w_2}$ and, using the
Morse theoretic approach, one shows that they are in fact connected
components. When $n=2$, even more components appear since, when $w_1=0$,
there is a reduction to the circle $\SO(2,\C) \subset
\mathrm{O}(2,\C)$ and this give rise to an integer invariant because
$\SO(2) =S^1$.

\begin{remark}
  These new invariants have been studied (and generalized) from the
  point of view surface group representations in the work of
  Guichard--Wienhard \cite{guichard-wienhard:2010}
\end{remark}

\begin{remark}
  Let $(W,Q)$ be the $\mathrm{O}(n,\C)$-bundle arising from a maximal
  $\Sp(2n,\R)$-Higgs bundle as above and define $\theta=(\beta\otimes
  1_{K^{n/2}})\circ Q\colon W \to W\otimes K^2$. Then $((W,Q),\theta)$
  is a $\GL(n,\R)$-Higgs bundle, except for the fact the twisting is
  by the square of the canonical bundle rather than the canonical
  bundle itself. This observation is the beginning of an interesting
  story known as the ``Cayley correspondence''; for more on this we
  refer to \cite{bradlow-garcia-prada-gothen:hss-higgs} and Rubio
  \cite{rubio:2012}. 
\end{remark}

\providecommand{\bysame}{\leavevmode\hbox to3em{\hrulefill}\thinspace}
\providecommand{\MR}{\relax\ifhmode\unskip\space\fi MR }
\providecommand{\MRhref}[2]{%
  \href{http://www.ams.org/mathscinet-getitem?mr=#1}{#2}
}
\providecommand{\href}[2]{#2}


\begin{thebibliography}{99}

\bibitem{bialynicki-birula:1973}
A.~Bia{\l}ynicki-Birula, \emph{Some theorems on actions of algebraic groups},
  Ann. of Math. (2) \textbf{98} (1973), 480--497. \MR{0366940 (51 \#3186)}

\bibitem{biswas-ramanan:1994}
I.~Biswas and S.~Ramanan, \emph{An infinitesimal study of the moduli of
  {H}itchin pairs}, J. London Math. Soc. (2) \textbf{49} (1994), 219--231.

\bibitem{bradlow-garcia-prada-gothen:2003}
S.~B. Bradlow, O.~Garc{\'\i}a-Prada, and P.~B. Gothen, \emph{Surface group
  representations and $\mathrm{U}(p,q)$-{H}iggs bundles}, J. Differential Geom.
  \textbf{64} (2003), 111--170.

\bibitem{bradlow-garcia-prada-gothen:hss-higgs}
\bysame, \emph{Maximal surface group representations in isometry groups of
  classical hermitian symmetric spaces}, Geometriae Dedicata \textbf{122}
  (2006), 185--213.

\bibitem{bradlow-garcia-prada-gothen:homotopy}
\bysame, \emph{Homotopy groups of moduli spaces of representations}, Topology
  \textbf{47} (2008), 203--224.

\bibitem{bradlow-garcia-gothen:2009}
\bysame, \emph{Deformations of maximal representations in
  $\mathrm{Sp}(4,\mathbb{R})$}, Q. J. Math. (2011), first published online June
  21, 2011. doi: 10.1093/qmath/har010.

\bibitem{bradlow-garcia-prada-mundet:2003}
S.~B. Bradlow, O.~Garc{\'{\i}}a-Prada, and I.~Mundet~i Riera, \emph{Relative
  {H}itchin-{K}obayashi correspondences for principal pairs}, Q. J. Math.
  \textbf{54} (2003), 171--208.

\bibitem{burger-iozzi-wienhard:2003}
M.~Burger, A.~Iozzi, and A.~Wienhard, \emph{Surface group representations with
  maximal {T}oledo invariant}, C. R. Math. Acad. Sci. Paris \textbf{336}
  (2003), no.~5, 387--390.

\bibitem{burger-iozzi-wienhard:2010b}
\bysame, \emph{Higher {T}eichm\"uller spaces: from $\mathrm{SL}(2,\mathbb{R})$
  to other {L}ie groups}, Handbook of {T}eichm\"uller Theory III, IRMA Lectures
  in Mathematics and Theoretical Physics, European Math.\ Soc., 2010, to
  appear.

\bibitem{burger-iozzi-wienhard:2010}
\bysame, \emph{Surface group representations with maximal {T}oledo invariant},
  Ann. of Math. (2) \textbf{172} (2010), no.~1, 517--566.

\bibitem{casimiro-florentino:2011}
A.~Casimiro and C.~Florentino, \emph{Stability of affine g-varieties and
  irreducibility in reductive groups}, Int.\ J. Math \textbf{23} (2012).

\bibitem{Chuang:2010ii}
W.-Y. Chuang, D.-E. Diaconescu, and G.~Pan, \emph{{Wallcrossing and Cohomology
  of The Moduli Space of Hitchin Pairs}}, Commun.Num.Theor.Phys. \textbf{5}
  (2011), 1--56.

\bibitem{corlette:1988}
K.~Corlette, \emph{Flat ${G}$-bundles with canonical metrics}, J. Differential
  Geom. \textbf{28} (1988), 361--382.

\bibitem{domic-toledo:1987}
A.~Domic and D.~Toledo, \emph{The {G}romov norm of the {K}aehler class of
  symmetric domains}, Math. Ann. \textbf{276} (1987), 425--432.

\bibitem{donaldson:1987}
S.~K. Donaldson, \emph{Twisted harmonic maps and the self-duality equations},
  Proc. London Math. Soc. (3) \textbf{55} (1987), 127--131.

\bibitem{dupont:1978}
J.~L. Dupont, \emph{Bounds for characteristic numbers of flat bundles},
  Springer LNM 763, 1978, pp.~109--119.

\bibitem{fock-goncharov:2006}
V.~V. Fock and A.~B. Goncharov, \emph{Moduli spaces of local systems and higher
  {T}eichmuller theory}, Publ. Math. Inst. Hautes \'Etudes Sci. \textbf{103}
  (2006), 1--211.

\bibitem{frankel:1959}
T.~Frankel, \emph{Fixed points and torsion on {K}\"{a}hler manifolds}, Ann. of
  Math. (2) \textbf{70} (1959), 1--8.

\bibitem{garcia-prada:2009}
O.~Garc{\'{\i}}a-Prada, \emph{Higgs bundles and surface group representations},
  Moduli spaces and vector bundles, London Math. Soc. Lecture Note Ser., vol.
  359, Cambridge Univ. Press, Cambridge, 2009, pp.~265--310.

\bibitem{garcia-gothen-mundet:2009a}
O.~Garc{\'{\i}}a-Prada, P.~B. Gothen, and I.~Mundet~i Riera, \emph{Higgs
  bundles and surface group representations in the real symplectic group},
  preprint, 2012, \texttt{arXiv:0809.0576v4 [math.AG]}.

\bibitem{garcia-gothen-mundet:2009b}
\bysame, \emph{The {H}itchin-{K}obayashi correspondence, {H}iggs pairs and
  surface group representations}, preprint, 2012, \texttt{arXiv:0909.4487v3
  [math.AG]}.

\bibitem{garcia-gothen-munoz:2004}
O.~Garc{\'\i}a-Prada, P.~B. Gothen, and V.~Mu{\~{n}}oz, \emph{Betti numbers of
  the moduli space of rank 3 parabolic {H}iggs bundles}, Mem. Amer. Math. Soc.
  \textbf{187} (2007), no.~879, viii+80.

\bibitem{garcia-heinloth-schmitt:2011}
O.~Garc{\'\i}a-Prada, J.~Heinloth, and A.~Schmitt, \emph{On the motives of
  moduli of chains and {H}iggs bundles}, \texttt{arXiv:1104.5558v1 [math.AG]},
  2011.

\bibitem{garcia-prada-mundet:2004}
O.~Garc{\'\i}a-Prada and I.~Mundet i~Riera, \emph{Representations of the
  fundamental group of a closed oriented surface in
  $\mathrm{Sp}(4,\mathbb{R})$}, Topology \textbf{43} (2004), 831--855.

\bibitem{goldman:1985}
W.~M. Goldman, \emph{Representations of fundamental groups of surfaces},
  Springer LNM 1167, 1985, pp.~95--117.

\bibitem{goldman:1988}
\bysame, \emph{Topological components of spaces of representations}, Invent.
  Math. \textbf{93} (1988), 557--607.

\bibitem{goldman:2010b}
W.~M. Goldman, \emph{Higgs bundles and geometric structures on surfaces}, The
  many facets of geometry, Oxford Univ. Press, Oxford, 2010, pp.~129--163.

\bibitem{goldman:2010}
W.~M. Goldman, \emph{Locally homogeneous geometric manifolds}, Proceedings of
  the International Congress of Mathematicians 2010 (R.~Bhatia, A.~Pal,
  G.~Rangarajan, V.~Srinivas, and M.~Vanninathan, eds.), World Scientific,
  2011.

\bibitem{goldman:1984}
William~M. Goldman, \emph{The symplectic nature of fundamental groups of
  surfaces}, Adv. in Math. \textbf{54} (1984), 200--225.

\bibitem{gothen:2001}
P.~B. Gothen, \emph{Components of spaces of representations and stable
  triples}, Topology \textbf{40} (2001), 823--850.

\bibitem{gothen-oliveira:2012}
P.~B. Gothen and A.~Oliveira, \emph{Rank two quadratic pairs and surface group
  representations}, Geometriae Dedicata (2012), first published online: 16
  March 2012. doi:10.1007/s10711-012-9709-1.

\bibitem{guichard-wienhard:2010}
O.~Guichard and A.~Wienhard, \emph{Topological invariants of {A}nosov
  representations}, J. Topol. \textbf{3} (2010), no.~3, 578--642.

\bibitem{hausel:2011}
T.~Hausel, \emph{Global topology of the {H}itchin system},
  \texttt{arXiv:1102.1717v2 [math.AG]}, 2011.

\bibitem{hernandez:1991}
L.~Hern{\'a}ndez, \emph{Maximal representations of surface groups in bounded
  symmetric domains}, Trans. Amer. Math. Soc. \textbf{324} (1991), 405--420.

\bibitem{hitchin:1987a}
N.~J. Hitchin, \emph{The self-duality equations on a {R}iemann surface}, Proc.
  London Math. Soc. (3) \textbf{55} (1987), 59--126.

\bibitem{hitchin:1987b}
\bysame, \emph{Stable bundles and integrable systems}, Duke Math. J.
  \textbf{54} (1987), 91--114.

\bibitem{hitchin:1992}
\bysame, \emph{{L}ie groups and {T}eichm\"{u}ller space}, Topology \textbf{31}
  (1992), 449--473.

\bibitem{mundet:2000}
I.~Mundet i~Riera, \emph{A {H}itchin-{K}obayashi correspondence for {K}\"ahler
  fibrations}, J. Reine Angew. Math. \textbf{528} (2000), 41--80.

\bibitem{kapustin-witten:2007}
A.~Kapustin and E.~Witten, \emph{Electric-magnetic duality and the geometric
  {L}anglands program}, Commun. Number Theory Phys. \textbf{1} (2007), 1--236.

\bibitem{kirwan:1984}
F.~Kirwan, \emph{Cohomology of quotients in symplectic and algebraic geometry},
  Mathematical Notes, vol.~31, Princeton University Press, Princeton, NJ, 1984.

\bibitem{labourie:1991}
F.~Labourie, \emph{Existence d'applications harmoniques tordues \`a valeurs
  dans les vari\'et\'es \`a courbure n\'egative}, Proc. Amer. Math. Soc.
  \textbf{111} (1991), no.~3, 877--882.

\bibitem{milnor:1957}
J.~W. Milnor, \emph{On the existence of a connection with curvature zero},
  Comment. Math. Helv. \textbf{32} (1958), 216--223.

\bibitem{narasimhan-seshadri:1965}
M.~S. Narasimhan and C.~S. Seshadri, \emph{Stable and unitary vector bundles on
  a compact {R}iemann surface}, Ann. Math. \textbf{82} (1965), 540--567.

\bibitem{rubio:2012}
Roberto~Rubio N\'u{\~n}ez, \emph{Higgs bundles and {H}ermitian symmetric
  spaces}, Ph.D. thesis, ICMAT -- Universidad Aut\'onoma de Madrid, 2012.

\bibitem{ramanathan:1975}
A.~Ramanathan, \emph{Stable principal bundles on a compact {R}iemann surface},
  Math. Ann. \textbf{213} (1975), 129--152.

\bibitem{schmitt:2008}
A.~Schmitt, \emph{Geometric invariant theory and decorated principal bundles},
  Z\"urich Lectures in Advanced Mathematics, European Mathematical Society,
  2008.

\bibitem{simpson:1988}
C.~T. Simpson, \emph{Constructing variations of {H}odge structure using
  {Y}ang-{M}ills theory and applications to uniformization}, J. Amer. Math.
  Soc. \textbf{1} (1988), 867--918.

\bibitem{simpson:1992}
\bysame, \emph{Higgs bundles and local systems}, Inst. Hautes {\'E}tudes Sci.
  Publ. Math. \textbf{75} (1992), 5--95.

\bibitem{toledo:1979}
D.~Toledo, \emph{Harmonic maps from surfaces to certain {K}aehler manifolds},
  Math. Scand. \textbf{45} (1979), 13--26.

\bibitem{toledo:1989}
\bysame, \emph{Representations of surface groups in complex hyperbolic space},
  J. Differential Geom. \textbf{29} (1989), 125--133.

\bibitem{turaev:1984}
V.~G. Turaev, \emph{A cocycle of the symplectic first {C}hern class and the
  {M}aslov index}, Funct. Anal. Appl. \textbf{18} (1984), 35--39.

\end{thebibliography}

\end{document}